\documentclass[preprint,11pt]{imsart}

\RequirePackage{amsthm,amsmath,amssymb,dsfont}
\RequirePackage[numbers]{natbib}
\RequirePackage[colorlinks,citecolor=blue,urlcolor=blue]{hyperref}

\usepackage{geometry}

\theoremstyle{plain}
\newtheorem{theorem}{Theorem}
\newtheorem{definition}{Definition}
\newtheorem{lemma}{Lemma}
\newtheorem{exm}{Example}
\newtheorem{assume}{Assumption}
\newtheorem{remark}{Remark}
\newtheorem{corollary}{Corollary}[theorem]
\newtheorem{proposition}{Proposition}[theorem]

\allowdisplaybreaks

\begin{document}

\begin{frontmatter}
\title{\bf On properties of fractional posterior in generalized reduced-rank regression}
\runtitle{On fractional posterior in generalized reduced-rank regression}
\begin{aug}
	\author{\fnms{The Tien}~\snm{Mai}\ead[label=e1]{the.t.mai@ntnu.no}\orcid{0000-0002-3514-9636}}
	\address{
		Department of Mathematical Sciences, 
		Norwegian University of Science and Technology,
		\\
		Trondheim 7034, Norway.
		\\
		\printead[presep={\ }]{e1}
	}
	\runauthor{T.T. Mai}
\end{aug}

\begin{abstract}
Reduced rank regression (RRR) is a widely employed model for investigating the linear association between multiple response variables and a set of predictors. While RRR has been extensively explored in various works, the focus has predominantly been on continuous response variables, overlooking other types of outcomes. This study shifts its attention to the Bayesian perspective of generalized linear models (GLM) within the RRR framework. In this work, we relax the requirement for the link function of the generalized linear model to be canonical. We examine the properties of fractional posteriors in GLM within the RRR context, where a fractional power of the likelihood is utilized. By employing a spectral scaled Student prior distribution, we establish consistency and concentration results for the fractional posterior. Our results highlight adaptability, as they do not necessitate prior knowledge of the rank of the parameter matrix. These results are in line with those found in frequentist literature. Additionally, an examination of model mis-specification is undertaken, underscoring the effectiveness of our approach in such scenarios.
\end{abstract}


\end{frontmatter}

\section{Introduction}
\label{sc_{ij}ntro}
The exploration of the relationship between multiple response variables and predictor sets has been extensively investigated in the literature. Reduced rank regression (RRR) has emerged as a prominent model for this purpose, leveraging a low-rank constraint to establish linear connections between response variables and predictors \citep{anderson1951estimating,izenman1975reduced,izenman2008modern,cook2018introduction,reinsel2023multivariate,giraud2021introduction}. This technique addresses the challenge of estimating coefficients in situations with large numbers of predictors or response variables. RRR effectively reduces the dimensionality of the problem, enabling accurate estimation even in high-dimensional scenarios \citep{bunea2011optimal,bunea2012joint}. Additionally, the reduced rank structure allows for an interpretation in terms of latent variables that explain the variation in response variables. Various estimation methods have been proposed for RRR, including both frequentist and Bayesian approaches, as in \cite{geweke1996bayesian,kleibergen2002priors,corander2004bayesian,schmidli2019bayesian,alquier2013bayesian,negahban2011estimation,rohde2011estimation,kleibergen2002priors,chen2013reduced,she2017robust,bing2019adaptive}. Furthermore, extensions of RRR, such as incorporating sparsity constraints, have also been investigated \citep{goh2017bayesian,chakraborty2020bayesian,yang2020fully}. This comprehensive exploration underscores the ongoing interest and effort dedicated to analyzing and modeling the relationships between multiple response variables and predictors.

However, a majority of these studies concentrate on continuous response data and utilize linear models. In practical scenarios, outcomes often vary in type, encompassing continuous measurements, binaries, categorical and counts data. Some recent works have begun focusing on multiple binary responses \citep{park2022low} or multiple count responses \citep{mishra2022negative}. Additionally, several efforts have been made to utilize the generalized linear model (GLM) framework, as in \cite{she2013reduced,luo2018leveraging,mishra2021generalized}. However, these aforementioned references predominantly adhere to frequentist principles, lacking a Bayesian perspective to address GLM in the RRR context. We endeavor to bridge this gap by introducing a Bayesian flavor to the analysis of GLM in the context of reduced rank regression. Moreover, in this work, we also relax the requirement for the link function of the generalized linear model to be canonical which is also a significant contribution.

Our emphasis is on investigating the theoretical properties of the fractional posterior within our framework, where a fractional power of the likelihood is employed, as elaborated in \cite{bhattacharya2016bayesian, alquier2020concentration}. This methodology is commonly known as generalized Bayesian approaches. Furthermore, the aforementioned references have demonstrated that fractional posteriors offer a more robust solution to tackling the issue of model misspecification. It is worth noting that generalized Bayesian inference has garnered increased interest in recent years, as evidenced by numerous studies, such as \cite{matsubara2022robust,hammer2023approximate,jewson2022general,yonekura2023adaptation, medina2022robustness,mai2017pseudo,grunwald2017inconsistency,bissiri2013general,yang2020alpha,lyddon2019general,syring2019calibrating,Knoblauch,mai2023reduced,hong2020model}. It is worth mentioning that fractional posteriors have been employed in the context of reduced-rank regression in \cite{chakraborty2020bayesian}, where the authors further incorporate sparsity considerations.

In this work, we establish consistency results for the fractional posterior in the context of generalized reduced rank regression, employing a low-rank promoting prior distribution on the matrix parameter. Our analysis yields consistency results in terms of the $\alpha$-R\'enyi divergence regarding the joint distribution of predictor and response variables. Subsequently, we extend our investigation to derive consistency results relative to the Hellinger metric and total variation distance. Additionally, we present consistency results for the matrix parameter directly, utilizing the Frobenius norm, contingent upon additional assumptions.

Furthermore, we explore concentration rates for the fractional posterior, considering the joint distribution of predictors and response variables. Our analysis extends to include concentration results concerning the $\alpha$-R\'enyi divergence, Hellinger metric, and total variation distance, as well as for the matrix parameter using the Frobenius norm. Remarkably, our findings indicate that these concentration rates adapt to the unknown rank of the true matrix parameter. We also establish consistency and concentration results for the mean estimator. Importantly, all these findings represent novel contributions within the scope of our knowledge. These results align with those found in frequentist literature. 

Moreover, we also investigate a scenario involving model mis-specification and demonstrate that both the fractional posterior and its mean estimator remain effective under this circumstance. The main result is presented through an oracle-type inequality for generalized Bayesian approaches in generalized reduced rank regression. In this work, we utilize a spectral scaled Student prior distribution, as it offers a computationally promising alternative via the gradient sampling method, Langevin Monte Carlo methods \cite{dalalyan2020exponential,mai2023reduced}. It should be noted that our results also apply to other types of priors based on low-rank factorization methods.

The rest of the paper is given as follows. In Section \ref{sc_model_method}, we introduce a GLM model for the reduced-rank regression problem, the prior and the fractional posterior are also given here. Section \ref{sc_main_results} presents  main results on the fractional posterior consistency and concentrations as well as results on estimations. The technical proofs are gathered in Section \ref{sc_proofss}. We conclude our work in Section \ref{sc_discusconclue}.

\section{Model and method}
\label{sc_model_method}
\subsection{Notations}
Let $\alpha\in(0,1)$ and $P,R$ be two probability measures. Let $\mu$ be any measure such that $P\ll \mu$ and $R\ll \mu$. The $\alpha$-R\'enyi divergence  between two probability distributions $P$ and $R$ is  defined by
$
D_{\alpha}(P,R)
 =
\frac{1}{\alpha-1} \log \int \left(\frac{{\rm d}P}{{\rm d}\mu}\right)^\alpha \left(\frac{{\rm d}R}{{\rm d}\mu}\right)^{1-\alpha} {\rm d}\mu  
,
$
and the Kullback-Leibler (KL) divergence is defined by
$
KL(P,R) 
= 
\int \log \left(\frac{{\rm d}P}{{\rm d}R} \right){\rm d}P$ if  $P \ll R
$, and $
+ \infty \text{ otherwise}.
$
The spectral norm, the Fobenius norm and the nuclear norm of a matrix $ A \in \mathbb{R}^{p\times q} $ will be respectively denoted by $ \|A\| , \|A\|_F  $ and $ \|A\|_* $. For $ d\geq 1 $, the identity matrix in $ \mathbb{R}^{d\times d} $ is denoted by $ \mathbf{I}_{d} $

\subsection{Model}

Let $ Y = [ y_1, \ldots, y_n ]^\top \in \mathbb{R}^{n\times q} $ be the response matrix consisting of $ n $ independent and identically distributed (i.i.d) observations from $ q $ response variables. More specifically, we consider the natural exponential family \citep{mccullagh1989generalized} that with a natural parameter $ \theta $ lying in $ \Theta \subset \mathbb{R} $  and a known scaling parameter $ a $, the probability density function of each $ y_{ij} $ takes the following form:
\begin{equation}
\label{eq_main_model}
p_{\theta_{ij}} (y_{ij} )
=
\exp\left\{ 
\frac{ y_{ij}\theta_{ij} -b(\theta_{ij}) }{a}
+
c(y_{ij} , a)
\right\}
.
\end{equation}
The function $ b(\cdot) $ is assumed to be twice-differentiable. In this model, one has that $ \mathbb{E} (Y_{ij}) = b'(\theta_{ij}) $ and $ Var(Y_{ij}) = ab''(\theta_{ij}) $, see \cite{mccullagh1989generalized}. Some of the most common distributions in this family include Bernoulli, Poisson and Gaussian (with known error variance) distributions. 

Let $ X = [ x_1, \ldots, x_n ]^\top \in \mathbb{R}^{n\times p} $ be the observed covariate/predictor matrix, where the number $ \, p\, $ of covariates can be much larger than the sample size $ n $. In this work, we assumed that the design matrix is fixed. We consider an increasing link function, denoted as $ h $, operating on a fixed covariate $ X_{i\cdot} \in \mathbb{R}^{ p} $ and a regression vector $ \beta_j \in \mathbb{R}^{ p} $, as in \cite{jeong2021posterior}. This function $ h $ satisfies the condition that its composition with $ b' $, denoted as $ (h\circ b') $, maps from the parameter space $ \Theta $ to the real numbers in a strictly increasing manner. This relationship is expressed as:
\begin{align}
\label{eq_linkGLM_noncanon}
(h\circ b') (\theta_{ij}) 
= 
X_{i\cdot}^\top \beta_j
,
\end{align}
where $ (h\circ b') (\cdot) = h[b'(\cdot)] $.
Here, we put $ B = [ \beta_1, \ldots, \beta_p ]^\top \in \mathbb{R}^{p\times q} $ as the matrix of regression coefficients for the predictors. For clarity, we will denote by $ B_0 $ the true matrix and use $ B $ as a notation for a generic parameter instead of $ \theta $ and put $ P_B $ be the corresponding distribution.

Put $ r^* : = {\rm rank} (B_0 ) $. Similar to the inquiries conducted in \cite{she2013reduced,luo2018leveraging} and \cite{park2022low}, we investigate scenarios wherein the matrix $ B_0 $ adheres to a low-rank structure ($ r^* \ll \min(p,q) $) or can be effectively approximated by a low-rank matrix. This assumption of low-rankness implies that the response matrix $ Y $ is predominantly influenced by the predictors through a restricted number of latent factors, see \cite{luo2018leveraging}.

In our context, the paper \cite{park2022low} has examined logistic regression as a specific case of the model described by equation \eqref{eq_main_model}. Additionally, the reference \cite{luo2018leveraging} has explored a broader interpretation of this model with a canonical link function, considering various exponential distributions for the $ q $ outcomes. Notably, these papers predominantly utilize frequentist methodologies. Therefore, we are pioneers in introducing a Bayesian perspective to this issue. Furthermore, in this work, we presents important and significant contributions  by considering a non-canonical link function. Notably, the formulation described in \eqref{eq_linkGLM_noncanon} does not necessitate $ h $ to be the inverse function of $ b' $, thus allowing for the inclusion of generalized linear models featuring non-canonical link functions. This aspect is crucial as our results apply to many models with non-canonical links, including probit regression, negative binomial regression, and gamma regression with the logarithmic link. Some detailed examples are given below.

\begin{exm}[Logistic regression]
	In this model, the responses are binaries and that $ b(\theta) = \log(1+e^\theta) , b'(\theta) = e^\theta/(1+e^\theta) $ and the link function is $ h(\theta) = \log[\theta/(1-\theta)] $.
\end{exm}

\begin{exm}[Probit regression]
	In this model, the responses are binaries and $ b(\theta) = \log(1+e^\theta) $. However, the link function takes the form $ h(\theta) = \Phi^{-1}(\theta) $ where $ \Phi $ is the distribution function of the standard normal distribution.
\end{exm}

\begin{exm}[Gamma regression with the log-link]
	In this case, the responses are positive values. For gamma regression with a known shape parameter, the natural parameter is the negative rate parameter. One has that $ b(\theta) = -\log(-\theta) $ and the link function is $ h(\theta) = \log(\theta) $.
\end{exm}

\begin{exm}[Poison regression]
	In this model, the responses are counts and $ b(\theta) = e^\theta $. The link function takes the form $ h(\theta) = \log(\theta) $.
\end{exm}

\begin{exm}[Negative binomial regression with the log-link]
	In this model, the responses are counts. With a given number of failures $ k $, negative binomial regression for the number of successes comes with $ b(\theta) = -k \log (1-e^\theta) $ and that  $ b'(\theta) = k e^\theta/ (1-e^\theta) $. The link function takes the form $ h(\theta) = \log(\theta) $.
\end{exm}

\subsection{Prior distribution and fractional posterior}

Let $\pi $ be a prior distribution for $ B $. The likelihood is denoted by
$
L_n(B) = \prod_{i=1}^n p_{B}(Y_{i})
.
$
In this work, we study the following fractional posterior, which serves as our \textit{ideal} estimator, is defined as
\begin{equation}
\label{eq_fractional_posterior}
\pi_{n,\alpha}(B )
\propto 
L_n^{\alpha}(B) \pi({\rm d}B),
\end{equation}
for some $ \alpha \in (0,1) $,
as in \cite{bhattacharya2016bayesian,alquier2020concentration}. The fractional posterior in \eqref{eq_fractional_posterior} yields the classical posterior distribution by taking $ \alpha =1 $. We define the mean estimator as
\begin{equation}
\label{eq_mean_estimator}
\hat{B} := 
\int B \pi_{n,\alpha}({\rm d}B)
.
\end{equation}

In this work, we have opted to use a spectral scaled Student distribution as a prior distribution. With a parameter $\tau>0$, it is given as
\begin{align}
\label{prior_scaled_Student}
\pi(B)
\propto
\det (\tau^2 \mathbf{I}_{p} + BB^\intercal )^{-(p+q+2)/2}.
\end{align}

\noindent This prior can induce low-rankness of matrices $ B $, as it can be verified that
$
\pi(B)
\propto
\prod_{j=1}^{p} (\tau^2 + s_j(B)^2 )^{- (p+q+2)/2 }
$, where $ s_j(B) $ denotes the $j^{th}$ largest singular value of $ B $. It means that this prior follows a scaled Student distribution evaluated at $ s_j(B) $ which induces approximately sparsity on the $s_j(B)$, see \cite{dalalyan2012sparse,dalalyan2012mirror}. Thus, under this prior distribution, most of the $s_j(B)$ are close to $0$ and that $B $ is approximately low-rank. This prior has been used before in different problem with matrix parameters as in \cite{yang2018fast,dalalyan2020exponential,mai2023bilinear,mai2023reduced}. Even though this prior distribution is not conjugate in our problem, it is advantageous to utilize the Langevin Monte Carlo, a sampling method that relies on gradients for implementation purposes.

\section{Main results}
\label{sc_main_results}
\subsection{Assumptions}
We adopt the following assumption regarding the second derivative $ b''(\cdot) $. These assumptions have been previously employed in the context of reduced rank regression, as outlined in \cite{luo2018leveraging}.
\begin{assume}
	\label{assume_boundon2ndorder}
	Assume that there exist a positive constant $  C_U < \infty $ such that the function $ b''(\cdot) $ satisfies the following conditions:
$
	\sup_{s\in\mathbb{R}}  b''(s) \leq C_U
	.
$
\end{assume}

\begin{assume}
	\label{assume_lowerbound2order}	
	Let's assume the existence of a  constant \( C_L >0\) such that \( C_L \leq C_U < \infty \), and the function \( b''(\cdot) \) fulfills the following criteria:
$
\inf_{s\in \Theta }  b''(s) \geq C_L
,
$
where $ \Theta $ is a closed (finite or infinite) interval in $ \mathbb{R} $.
\end{assume}

For Gaussian distributions with $ b(\theta) = \theta^2/2 $, the second derivative $ b''(\theta) = 1 $, resulting in both upper and lower  constraints, $ C_U = C_L = 1 $ for all $ \theta $. 
When considering the binomial distribution, $ b''(\theta) = e^\theta / (1+e^\theta)^2 $, leading to $ C_U = 1/4 $.
However, Assumption \ref{assume_lowerbound2order} is equivalent to constraining the range of $ \theta $ within $ \Theta = \{\theta: |\theta|\leq C_0 \} $, which in turn yields $ C_L = e^{C_0} / (1+e^{C_0})^2 $. Remarkably, Assumption \ref{assume_boundon2ndorder} aims to preclude scenarios where the variance $ \text{Var}(Y) $ is unbounded, while Assumption \ref{assume_lowerbound2order} is designed to avoid situations where the variance is exceedingly small.

\subsection{Consistency results}
\label{sc_consistency}
Initially, we present consistency results concerning our fractional posterior. Subsequently, concentration rates are expounded upon in Section \ref{sc_concetration}. To the best of our knowledge, these results represent entirely novel contributions.

\begin{theorem}
	\label{thm_result_expectation}
	For any $\alpha\in(0,1)$, under Assumption \ref{assume_boundon2ndorder}, and for $ \tau^2 = 2a/(qp \|X \|_{F}^2) $,
	we have that
	\begin{equation*}
	\mathbb{E} \left[ \int D_{\alpha}(P_{B},P_{B_0}) \pi_{n,\alpha}({\rm d} B ) \right]
	\leq 
	\frac{1+\alpha}{1-\alpha}\varepsilon_n
	.
	\end{equation*}
	where
$$
	\varepsilon_n
	=
	\frac{ C_U 2 r^* (q+p+2) \log \left(  1+ \frac{\| X \|_F \| B_0 \|_F  \sqrt{qp}}{\sqrt{4a r^*}} \right) 
	}{nq}
	.
$$
\end{theorem}

We remind that all technical proofs are deferred to Section \ref{sc_proofss}. It is noteworthy that in the rate $ \varepsilon_n $, the condition $ r^* = {\rm rank} (B_0 ) \neq 0 $ is not required. If $ r^* = 0 $, then $ B_0 = 0 $, and we interpret $ 0\log(1+0/0) $ as $ 0 $.

\begin{remark}
	It is important to note that our results are formulated without prior knowledge of $ r^* $, the rank of the true underlying parameter matrix.
The implications drawn from Theorem \ref{thm_result_expectation} suggest that the consistency of the model adjusts to the unknown rank of the true parameter matrix $ B_0 $. This finding present a novel establishment of consistency for the fractional posterior in generalized reduced rank regression. The rate $ \varepsilon_n $ scales with the order of $ r^*(q+p)/n $, up to a logarithmic factor. The logarithmic term can be further simplified to $ \log ( \small 1+ \| X \|_F \| B_0 \| \sqrt{pq} ) $, utilizing the inequality $ \| B_0 \|_F  \leq \| B_0 \| \sqrt{r^*} $.
\end{remark}

\begin{remark}
 As specific instances, we deduce the following noteworthy corollary on the Hellinger distance and the total variation distance by leveraging results from \cite{van2014renyi}. For example, with $\alpha\in(0,1]$,
 $
 (\alpha/2) d^{2}_{TV}(P,R) \leq D_{\alpha}(P,R) 
 ,
 $
 $d_{TV}$ being the total variation distance. 	And
 $
 H^2(P,R) 
 =
 2[1-\exp(-(1/2)D_{1/2}(P,R))] 
 \leq
 D_{1/2}(P,R)
 $
 where $ H^2(\cdot,\cdot) $ is the squared Hellinger distance.
 The $\alpha$-R\'enyi divergences are all equivalent \citep{van2014renyi}, i.e. for $0<\alpha<1$ and $\alpha \leq \beta $,
 $
 \frac{\alpha}{\beta}\frac{1-\beta}{1-\alpha} D_{\beta} \leq D_\alpha \leq D_\beta 
 .
 $
 
\end{remark}

For sake of simplicity, we put
$
c_\alpha 
=
\frac{2(\alpha+1)}{1-\alpha}, \alpha \in [0.5,1) ,
$ and $
c_\alpha 
= \frac{2(\alpha+1)}{\alpha}, \alpha \in (0, 0.5)
$. We immediately obtain the following result.

\begin{corollary}
	\label{cor_concentration_Hellinger}
A special case of Theorem \ref{thm_result_expectation} yields a concentration result in terms of the classical Hellinger distance
	\begin{equation*}
\mathbb{E} \left[ \int H^2(P_{B},P_{B_0}) \pi_{n,\alpha}({\rm d} B ) \right]
\leq 
c_\alpha \varepsilon_n
.
\end{equation*}
	And, 
	\begin{equation*}
	\mathbb{E} \left[ \int  d^{2}_{TV}(P_{B},P_{B_0}) \pi_{n,\alpha}({\rm d} B ) \right]
	\leq 
	\frac{2(1+\alpha)}{\alpha(1-\alpha)}
	 \varepsilon_n
	,
	\end{equation*}
	with $d_{TV}$ being the total variation distance.
\end{corollary}

Although the above outcomes ensure the consistency of the fractional posterior concerning measures such as $\alpha$-R\'enyi divergence, Hellinger distance, or total variation distance between the densities, they do not address the parameter matrix directly. To address this gap, we proceed to furnish results regarding the matrix parameter estimation within certain specified distances. In order to do that, it becomes necessary to impose further constraints, as outlined in Assumption \ref{assume_lowerbound2order}.

\begin{proposition}
	\label{propos_estimation}
Under the assumptions outlined in Theorem \ref{thm_result_expectation} and with the additional requirement that Assumption \ref{assume_lowerbound2order} is satisfied, one has that
	\begin{equation*}
\mathbb{E} \left[ \int  \frac{1}{nq}
\| X^\top(B - B_0)\|_F^2 \pi_{n,\alpha}({\rm d} B ) \right]
\leq 
\frac{ 2a (1+\alpha)}{ C_L (1- \alpha) } 
 \varepsilon_n
,
\end{equation*}		
and consequently
	\begin{equation*}
\mathbb{E} \left[    \frac{1}{nq}
\| X^\top(\hat{B} - B_0)\|_F^2  \right]
\leq 
\frac{ 2a (1+\alpha)}{ C_L (1- \alpha) }  \varepsilon_n
.
\end{equation*}	
\end{proposition}

To derive results for the mean estimator, as defined in \eqref{eq_mean_estimator}, we utilize Jensen's inequality with a convex function. Similarly, results for the median of the fractional posterior can be obtained using Jensen's inequality, as explained  in \cite{merkle2005jensen}.

\subsection{Result in the misspecified case}
\label{sc_mispecifiedcase}
In this section, we show that our results in Section \ref{sc_consistency} can be extended to the misspecified setting.

\begin{assume}
	\label{assume_bound1storder}
	Assume that there exist a positive constant $  U_1 < \infty $ such that the function $ b'(\cdot) $ satisfies the following conditions:
$
	\sup_{s\in\mathbb{R}}  b'(s) \leq U_1
	.
$
\end{assume}
This assumption is used to obtain the results for the case of model mis-specified. For Gaussian distributions with $ b(\theta) = \theta^2/2 $, Assumption \ref{assume_bound1storder} means that the parameter is assumed to be upper bounded by a constant.

Assume that the true data generating distribution is now parametrized by $ B_0\notin\mathcal{B} $, where $ \mathcal{B} $ is the parameter space considered in Section \ref{sc_consistency}. Let's define $P_{B_0}$ as the true distribution. 	Put $ \bar{B} = \arg\min_{ B\in \mathcal{B}} KL ( P_{B_0},P_{B} ) $.

\begin{theorem}
	\label{theorem_misspecified}
	For any $\alpha\in(0,1)$, let assume that Assumption \ref{assume_bound1storder} holds,  with $\tau = a/( 2 \sqrt{ n q }\sqrt{pq} \|X \|_{F}) $. Then, 
	\begin{equation*}
	\mathbb{E} \left[ \int D_{\alpha}(P_{B},P_{B_0}) \pi_{n,\alpha}({\rm d}B )\right]
	\leq 
	\frac{\alpha}{1-\alpha} \min_{ B \in \mathcal{B} } 
	KL (P_{B_0},P_{B})
	+ \frac{1+\alpha}{1-\alpha} r_n
	,
	\end{equation*}
	where 
	$$ 
	r_n = 	
\frac{	2 U_1 {\rm rank} (\bar{B}) (q+p+2) \log \left( 1+  \frac{\| X \|_F \| \bar{B}\|_F  2 \sqrt{ n q }\sqrt{pq} }{ a\sqrt{ 2 {\rm rank}(\bar{B})}} \right) 
}{ nq}
	.
	$$
\end{theorem}

In the case of a well-specified model, i.e., when $ B_0 = \arg\min_{ B\in B} KL ( P_{B_0},P_{B} ) $, we retrieve Theorem~\ref{thm_result_expectation}. Otherwise, this result manifests as an oracle inequality. Although the differing risk measures on both sides, this observation remains valuable, particularly when $ KL ( P_{B_0},P_{B} ) $ is minimal. 

Nonetheless, under additional assumptions, we can further derive an oracle inequality result with $\ell_2$ distance on both sides. The result is as follows.

\begin{corollary}
	\label{cor_estimation_miss}
	Under the same assumptions as in Theorem~\ref{theorem_misspecified} and with the additional assumption that Assumption \ref{assume_boundon2ndorder} and  Assumption \ref{assume_lowerbound2order} hold. For $ 0 < \alpha < 1 $, then
	\begin{multline*}
\mathbb{E} \left[ \int 	 \frac{1}{nq} 
\| X^\top (B-B_0) \|_F^2 \pi_{n,\alpha}({\rm d}B )\right]
\leq 
\inf_{0\leq r \leq pq}
\inf_{
	\begin{tiny}
	\begin{array}{c}
	\bar{B}\in\mathbb{R}^{p\times q}
	\\
	{\rm rank}(\bar{B}) \leq r
	\end{array}
	\end{tiny}
} 
\Biggl\{
\frac{C_U}{nqC_L} \frac{\alpha}{1-\alpha}
\| X^\top (\bar{B}-B_0) \|_F^2
\\
+ 
\frac{4a(1+\alpha)}{ C_L (1-\alpha)}
\frac{ U_1 r (q+p+2) \log 
\left( 1+  \frac{\| X \|_F \| \bar{B}\|_F  2 \sqrt{ n q }\sqrt{pq} }{ a\sqrt{ 2 r } } \right) 
}{ nq}
\Biggr\}
,
\end{multline*}
consequently, using Jensen's inequality,
	\begin{multline*}
\mathbb{E} 	  
\| X^\top (\hat{B} - B_0) \|_F^2 
\leq 
\inf_{0\leq r \leq pq}
\inf_{
	\begin{tiny}
	\begin{array}{c}
	\bar{B}\in\mathbb{R}^{p\times q}
	\\
	{\rm rank}(\bar{B}) \leq r
	\end{array}
	\end{tiny}
} 
\Biggl\{
\frac{C_U}{nqC_L} \frac{\alpha}{1-\alpha}
\| X^\top (\bar{B}-B_0) \|_F^2
\\
+ 
\frac{4a(1+\alpha)}{ C_L (1-\alpha)}
\frac{ U_1 r (q+p+2) \log 
	\left( 1+  \frac{\| X \|_F \| \bar{B}\|_F  2 \sqrt{ n q }\sqrt{pq} }{ a\sqrt{ 2 r } } \right) 
}{ nq}
\Biggr\}
.
\end{multline*}

\end{corollary}

Up to our knowledge, these results are  novel.

\subsection{Concentration results}
\label{sc_concetration}

We will now present stronger results concerning the concentration rates of the fractional posterior. In contrast to the consistency results discussed earlier, obtaining concentration results for our fractional posterior necessitates verifying additional conditions. In additional to checking the conditions,
$
\frac1n KL(\rho_n,\pi)
\leq
\varepsilon_n
,
$
and
$
\int KL(P_{\theta_0},P_{\theta})
\rho_n(d\theta)
\leq 
\varepsilon_n
,
$
as in the proof of Theorem \ref{thm_result_expectation}, we need to further verify that 
$
\int
\mathbb{E}\left[\log^2
\left(
\frac{p_{\theta}(Y)}{p_{\theta_0}(Y)}
\right)
\right] 
\rho_n(d\theta)
\leq
\varepsilon_n
.
$
This is detailed in Theorem 2.4 from \cite{alquier2020concentration}.

The concentration results provided by Theorem \ref{thm_result_probabily} for the fractional posterior in generalized reduced rank regression represent, to our knowledge, novel contributions to this problem.

\begin{theorem}
	\label{thm_result_probabily}
	For any $\alpha\in(0,1)$, under Assumption \ref{assume_boundon2ndorder}, and for $\tau^2 = a/(qp \|X \|_{F}^2)$,
	we have that
	\begin{equation*}
	\mathbb{P} \left[ \int D_{\alpha}(P_{B},P_{B_0}) \pi_{n,\alpha}({\rm d} B ) 
	\leq 
	2	\frac{1+\alpha}{1-\alpha}\varepsilon_n
	\right]
	\geq 
	1 -	\frac{2}{n\varepsilon_n}
	.
	\end{equation*}
	where
$$
	\varepsilon_n
	=
\frac{C_U^2}{4nq}
\vee
\frac{ 	2C_U r^* (q+p+2) \log \left( 1+ \frac{\| X \|_F \| B_0 \|_F  \sqrt{qp}}{\sqrt{ 2ar^*}}  \right) 
}{nq}
	,
$$
	with the convention $0\log(1+0/0)=0$.
\end{theorem}

The proof of this theorem is detailed in Section \ref{sc_proofss}. The proof methodology employed in this theorem follows a general approach outlined in \cite{alquier2020concentration} and \cite{bhattacharya2016bayesian}, which provide a comprehensive framework for analyzing fractional posteriors.

Similar to the results presented in Corollary \ref{cor_concentration_Hellinger}, we can derive concentration properties for the fractional posterior using both the Hellinger distance and the total variation distance.

\begin{corollary}
	\label{cor_concentration_Hellinger2}
In specific instances, Theorem \ref{thm_result_probabily} gives rise to concentration results expressed in terms of the classical Hellinger distance and the total variation, as outlined below:
	\begin{equation*}
	\mathbb{P} \left[ \int H^2(P_{B},P_{B_0}) \pi_{n,\alpha}({\rm d} B ) 
	\leq 
	c_\alpha \varepsilon_n
		\right]
	\geq 
	1 -	\frac{2}{n\varepsilon_n}
	,
	\end{equation*}
	\begin{equation*}
	\mathbb{P} \left[ \int  d^{2}_{TV}(P_{B},P_{B_0}) \pi_{n,\alpha}({\rm d} B ) 
	\leq 
	\frac{2(1+\alpha)}{\alpha(1-\alpha)}
	\varepsilon_n
		\right]
	\geq 
	1 -	\frac{2}{n\varepsilon_n}
	,
	\end{equation*}
	with $d_{TV}$ being the total variation distance.
\end{corollary}

In order to derive results related to prediction error \( \| X^\top( \hat{B} - B_0)\|_F^2 \), it is crucial to incorporate additional Assumption \ref{assume_lowerbound2order}.

\begin{proposition}
	\label{propos_estimation2}
Under the assumptions specified in Theorem \ref{thm_result_probabily}, and with the additional condition that Assumption \ref{assume_lowerbound2order} is met, we establish that
	\begin{equation*}
	\mathbb{P} \left[ \int  
		\frac{1}{nq}
	\| X^\top(B - B_0)\|_F^2 \pi_{n,\alpha}({\rm d} B ) 
	\leq 
	\frac{ 4a (1+\alpha)}{ C_L (1-\alpha) }  \varepsilon_n
		\right]
	\geq 
	1 -	\frac{2}{n\varepsilon_n}
	,
	\end{equation*}	
	and consequently,
	\begin{equation*}
\mathbb{P} \left[  
	\frac{1}{nq} 
\| X^\top( \hat{B} - B_0)\|_F^2
\leq 
\frac{ 4a (1+\alpha)}{ C_L (1-\alpha) }  \varepsilon_n
\right]
\geq 
1 -	\frac{2}{n\varepsilon_n}
.
\end{equation*}			
	
\end{proposition}

\begin{remark}
The minimax error rate for the linear model with a Gaussian distribution and known error variance, evaluated under the loss \( \| X^\top( \hat{B} - B_0)\|_F^2 \), was established in \cite{bunea2011optimal} to be of the order \( r^*(p+q) \). Our obtained rate aligns with this rate, up to a logarithmic factor, in the traditional reduced rank regression model. 	In the context of using a reduced-rank logistic regression model for multiple binary responses, our result in squared Frobenius norm is similar to those from \cite{park2022low}. For the generalized reduced rank regression, as discussed in \cite{luo2018leveraging}, our results from Proposition \ref{propos_estimation} and Proposition \ref{propos_estimation2} yield similar results, albeit with distinct logarithmic terms.
\end{remark}

In addition to achieving the optimal prediction rate, we are prepared to derive estimation rates but first need to introduce additional notation. Let's define:
\[
\kappa := \min_{ B \neq 0 } \frac{\| X^\top B\|_F}{\sqrt{n} \| B\|_F}
\]
Note that \( \kappa \) represents the restricted eigenvalue constant \citep{bickel2009simultaneous}. Similar notation has been employed for logistic reduced rank regression in \cite{park2022low}. With this notation, we readily obtain the following result.
\begin{corollary}
	\label{coro_esti_2}
	Under the assumptions specified in Proposition \ref{propos_estimation2}, and assuming that $ \kappa $ is strictly positive, then
	\begin{equation*}
	\mathbb{P} \left[ \int  
	\| B - B_0\|_F^2 \pi_{n,\alpha}({\rm d} B ) 
	\leq 
	\frac{ 4a (1+\alpha)}{ \kappa^2 C_L (1-\alpha) }  \varepsilon_n'
	\right]
	\geq 
	1 -	\frac{2}{n\varepsilon_n}
	,
	\end{equation*}	
	and consequently, using Jensen's inequality,
	\begin{equation*}
	\mathbb{P} \left[  
	\| \hat{B} - B_0 \|_F^2
	\leq 
	\frac{ 4a (1+\alpha)}{ \kappa^2 C_L (1-\alpha) }  \varepsilon_n'
	\right]
	\geq 
	1 -	\frac{2}{n\varepsilon_n}
	,
	\end{equation*}	
where	
	$$
	\varepsilon_n'
	=
	\frac{C_U}{an}
	\vee
	\frac{C_U^2}{4a^2n}
	\vee
	\frac{ 	2 r^* (q+p+2) \log \left( 1+ \frac{\| X \|_F \| B_0 \|_F  \sqrt{qp}}{\sqrt{ 2r^*}}  \right) 
	}{n}
	,
	$$		
	
\end{corollary}

\begin{remark}
The rate for \( \| \hat{B} - B_0 \|_F^2 \) of order \( r^*\max(p,q)/n \) is also akin to those found in \cite{park2022low} for a logistic model.
\end{remark}

\section{Proofs}
\label{sc_proofss}

\subsection{Proofs for Section \ref{sc_consistency}}

\begin{proof}[{\bf Proof of Theorem~\ref{thm_result_expectation}}]
We can verify the hypotheses regarding the Kullback-Leibler divergence between the likelihood terms, as required in Theorem 2.6 of \cite{alquier2020concentration}, for the expression of $\rho_n$ provided in equation \eqref{eq_priorspecific_anark}.

	First, from Lemma \ref{lema_KLupperbound}, we have that 
	\begin{align*}
KL(P_{\theta_0},P_{\theta} )
& \leq
\frac{C_U}{2a}\frac{1}{nq} \| \theta_0 -\theta \|_F^2
\\
& =
\frac{C_U}{2a}\frac{1}{nq} \| X^\top (B-B_0) \|_F^2
\\
& \leq
\frac{C_U}{2a}\frac{1}{nq} \| X\|_F^2 \| B-B_0 \|_F^2
.
\end{align*}
When integrating with respect to $\rho_n$, we have that
	\begin{align}
	\int KL(P_{\theta_0},P_{\theta})
	\rho_n(d\theta)
&	\leq 
\frac{C_U}{2a}\frac{1}{nq} \| X\|_F^2 \int \| B-B_0 \|_F^2	\rho_n(d B)
\nonumber
\\
& = 
\frac{C_U}{2a}\frac{1}{nq} \| X\|_F^2 \int \| B-B_0 \|_F^2	\pi(B - B_0) dB
\nonumber
\\
&	\leq 
\frac{C_U}{2a}\frac{1}{nq} \| X\|_F^2 qp \tau^2
\label{eq_bound_kl1}
	,
	\end{align}
where we have used Lemma \ref{lemma:arnak:1} in the last inequality and using the change of variable $ M = B-B_0 $. 

From Lemma \ref{lemma:arnak:2}, we have that 
	\begin{align*}
	\frac{1}{nq} KL(\rho_n,\pi)
	\leq
	\frac{	2r^* (q+p+2) \log \left( 1+ \frac{\| B_0 \|_F}{\tau \sqrt{2 r^*}} \right) 
	}{nq}
	.
	\end{align*}
For $\tau^2 = 2a/(qp \|X \|_{F}^2)$, it leads to that
	\begin{align*}
\frac{1}{nq} KL(\rho_n,\pi)
\leq
\frac{	2 r^* (q+p+2) \log \left( 1+  \frac{\| X \|_F \| B_0 \|_F  \sqrt{qp}}{\sqrt{ 4a r^* }} \right) 
}{nq}
,
\end{align*}
and
	\begin{align*}
\int KL(P_{\theta_0},P_{\theta})
\rho_n(d\theta)
	\leq 
 \frac{C_U}{nq}
.
\end{align*}
	
Therefore, we can apply Theorem 2.6 in \cite{alquier2020concentration} with $\rho_n$ provided in  \eqref{eq_priorspecific_anark} and
	\[
	\varepsilon_n
	=
	\frac{ C_U	2 r^* (q+p+2) \log \left( 1+ \frac{\| X \|_F \| B_0 \|_F  \sqrt{qp}}{\sqrt{ 4a r^*}}  \right) 
	  }{nq}
	.
	\]
	The proof is completed.
	
\end{proof}

\begin{proof}[\bf Proof of Proposition \ref{propos_estimation}]
	From Theorem \ref{thm_result_expectation}, one has that
	\begin{equation*}
	\mathbb{E} \left[ \int D_{\alpha}(P_{B},P_{B_0}) \pi_{n,\alpha}({\rm d} B ) \right]
	\leq 
	\frac{1+\alpha}{1-\alpha}\varepsilon_n
	,
	\end{equation*}
	and from Lemma \ref{lema_lowerboundreymi}, where $ \theta = X^\top B , \theta_0 = X^\top B_0 $, one has that
	\begin{align*}
	D_{\alpha}(P_{B},P_{B_0})
	\geq
	\frac{C_L }{2a }  \frac{1}{nq}
	\|X^\top B -X^\top B_0 \|_F^2
	.
	\end{align*}
	Thus, we obtain that
	\begin{equation*}
\frac{1}{nq}	\mathbb{E} \left[ \int  
	\| X^\top(B - B_0)\|_F^2 \pi_{n,\alpha}({\rm d} B ) \right]
	\leq 
	\frac{2a}{C_L}
	\frac{1+\alpha}{1- \alpha }
	\varepsilon_n
	.
	\end{equation*}	
By applying Jensen's inequality with a convex function, one has that
	\begin{equation*}
\| X^\top(\hat{B} - B_0)\|_F^2 
\leq
\int 
\| X^\top(B - B_0)\|_F^2 
\pi_{n,\alpha}({\rm d}M) 
,
\end{equation*}
consequently,
	\begin{equation*}
\frac{1}{nq}\mathbb{E} \left[   
\| X^\top(\hat{B} - B_0)\|_F^2 \right]
\leq 
	\frac{2a}{C_L}
\frac{1+\alpha}{ 1 - \alpha }
\varepsilon_n
.
\end{equation*}			
	The proof is completed.	
	
\end{proof}

\subsection{Proof of Section \ref{sc_mispecifiedcase}}

\begin{proof}[\bf Proof of Theorem~\ref{theorem_misspecified}]
	
	Put $ \bar{B} = \arg\min_{ B\in \mathcal{B}} KL ( P_{B_0},P_{B} ) $. We will verify the conditions outlined in Theorem 2.7 of \cite{alquier2020concentration} in order to apply it.
	
	First, we have from Lemma \ref{lema_KLup_missped} that
	\begin{align*}
	\mathbb{E}_{ B_0}\left[\log\frac{ p_{ \bar{B}} }{ p_{ B}}(Y)\right]
	& \leq
	\frac{2U_1}{a \sqrt{nq}} 
	\| X^\top\bar{B} - X^\top B \|_F
	\\
	& \leq
	\frac{2U_1}{a \sqrt{nq}} 
	\| X\|_F \| B- \bar{B} \|_F
	.	
	\end{align*}
	
	Consider the following distribution as a translation of the prior $ \pi $,
$
	\rho^*_{n} ( B) 
	\propto 
	\pi ( B -\bar{B})
	.
$
	Now,
	\begin{align*}
	\int \mathbb{E}_{ B_0}
	\left[\log\frac{ p_{ \bar{B}}}{ p_{ B}} (Y) \right]
	\rho^*_n({\rm d} B) 
	&	\leq 
	\frac{2U_1}{a \sqrt{nq}} 
	\| X\|_F \int \| B- \bar{B} \|_F
	\pi ( B -\bar{B}) {\rm d} B	
	\\
	&	\leq
	\frac{2U_1}{a \sqrt{nq}} 
	 \| X\|_F 
	\left(
	\int \| B- \bar{B} \|_F^2
	\pi ( B -\bar{B}) {\rm d} B	
	\right)^{1/2}	
	.
	\end{align*}
	where we have Holder's inequality in the last inequality above.
	Applying Lemma \ref{lemma:arnak:1} with a change of variables, one gets that
	\begin{align*}
	\int \mathbb{E}_{ B_0}
	\left[\log\frac{ p_{ \bar{B}}}{ p_{ B}} (Y) \right]
	\rho^*_n({\rm d} B) 
	\leq 
	\frac{2U_1}{a \sqrt{nq}} 
	\| X\|_F \sqrt{ q p \tau^2}
	.
	\end{align*}	
	From Lemma \ref{lemma:arnak:2}, we have that 
	\begin{align*}
	\frac{1}{nq} KL(\rho^*_n,\pi)
	\leq
	\frac{	2 {\rm rank} (\bar{B} ) (q+p+2) \log \left( 1+ \frac{\| \bar{B} \|_F}{\tau \sqrt{2{\rm rank} (\bar{B} )}} \right) 
	}{nq}
	.
	\end{align*}
	For $\tau = a/( 2 \sqrt{ n q }\sqrt{pq} \|X \|_{F}) $, it leads to that
	\begin{align*}
	\int \mathbb{E}_{ B_0}
	\left[\log\frac{ p_{ \bar{B}}}{ p_{ B}} (Y) \right]
	\rho^*_n({\rm d} B) 
	\leq 
	\frac{U_1}{nq}
	,
	\end{align*}
	and
	\begin{align*}
	\frac{1}{nq} KL(\rho^*_n,\pi)
	\leq
	\frac{	2 {\rm rank} (\bar{B}) (q+p+2) \log \left( 1+  \frac{\| X \|_F \| \bar{B}\|_F  2 \sqrt{ n q }\sqrt{pq} }{ a\sqrt{ 2 {\rm rank}(\bar{B})}} \right) 
	}{nq}
	.
	\end{align*}
	Consequently, we can apply Theorem 2.7 in \cite{alquier2020concentration} with $ \rho_n := \rho_n^* $ and with
	\[
	r_n
	=
	\frac{	2 U_1 {\rm rank} (\bar{B}) (q+p+2) \log \left( 1+  \frac{\| X \|_F \| \bar{B}\|_F  2 \sqrt{ n q }\sqrt{pq} }{ a\sqrt{ 2 {\rm rank}(\bar{B})}} \right) 
	}{ nq}
	.
	\]
	The proof is completed.
	
\end{proof}

\begin{proof}[\bf Proof of Corollary~\ref{cor_estimation_miss}]
	First, from Lemma \ref{lema_KLupperbound}, we have that 
	\begin{align*}
	KL(P_{B_0},P_{B} )
	\leq
	\frac{C_U}{2a} \frac{1}{nq} 
	\| X^\top (B-B_0) \|_F^2
	.
	\end{align*}
	Now, from Lemma \ref{lema_lowerboundreymi}, we have that
	\begin{align*}
	D_{\alpha}(P_{B},P_{B_0})
	\geq
	\frac{C_L }{2a }  \frac{1}{nq}  
	\| X^\top (B-B_0) \|_F^2
	.
	\end{align*}
	Plug in these bounds into Theorem \ref{theorem_misspecified}, it yields that
	\begin{equation*}
	\mathbb{E} \left[ \int 	
	\frac{C_L }{2a  }  
	\| X^\top (B-B_0) \|_F^2 \pi_{n,\alpha}({\rm d}B )\right]
	\leq 
	\frac{\alpha}{1-\alpha} \min_{ B \in B } 
	\frac{C_U}{2a} \| X^\top (B-B_0) \|_F^2
	+ \frac{1+\alpha}{1-\alpha} nq r_n
	,
	\end{equation*}
	consequently,
	\begin{equation*}
	\mathbb{E} \left[ \int 	  
	\| X^\top (B-B_0) \|_F^2 \pi_{n,\alpha}({\rm d}B )\right]
	\leq 
	\frac{C_U}{C_L}	\frac{\alpha}{1-\alpha}
	\min_{ B \in B } 
	\| X^\top (B-B_0) \|_F^2
	+ 
	\frac{2a(1+\alpha)}{ C_L(1-\alpha)} nq r_n
	.
	\end{equation*}
	The proof is completed.
	
\end{proof}

\subsection{Proofs for Section \ref{sc_concetration}}

\begin{proof}[{\bf Proof of Theorem \ref{thm_result_probabily}}]
	We will apply Theorem 2.4 in \cite{alquier2020concentration}. In order to do so, we need to verify the hypotheses regarding the Kullback-Leibler divergence between the likelihood terms as required for $\rho_n$ provided in equation \eqref{eq_priorspecific_anark}.
	
	From \eqref{eq_bound_kl1}, we have that
	\begin{align*}
	\int KL (P_{B_0},P_{B})
	\rho_n(dB)
	\leq 
	\frac{C_U}{2a} \frac{1}{nq}  \| X\|_F^2 qp \tau^2
	.
	\end{align*}
	Now,  from Lemma \ref{lm_bound_varkl}, where $ \theta = X^\top B , \theta_0 = X^\top B_0 $, one has that 
	\begin{align*}
	\mathbb{E}\left[\log
	\left(
	\frac{p_{B}}{p_{B_0} }(Y)
	\right)^2
	\right] 
	& \leq
	\frac{C_U}{a nq } \| X^\top (B-B_0) \|_F^2
	+
	\frac{C_U^2}{4a^2 nq } \| X^\top (B-B_0) \|_F^4
	\\
	& \leq
	\frac{C_U}{a nq } \| X\|_F^2 \| B-B_0 \|_F^2
	+
	\frac{C_U^2}{4a^2 nq } \| X\|_F^4 \|B-B_0 \|_F^4
	.
	\end{align*}
	When integrating with respect to $\rho_n$, we have 
	\begin{align*}
	&	\int
	\mathbb{E}\left[\log
	\left(
	\frac{p_{B}}{p_{B_0}}(Y)
	\right)^2
	\right] 
	\rho_n(dB)
	\\
	&	\leq 
	\frac{C_U}{a nq } \| X\|_F^2 \int \| B-B_0 \|_F^2	\rho_n(d B) 
	+
	\frac{C_U^2}{4a^2 nq } \| X\|_F^4 \int \|(B-B_0) \|_F^4 	\rho_n(d B)
	\\
	& = 
	\frac{C_U}{a nq } \| X\|_F^2 \int \| B-B_0 \|_F^2	\pi(B - B_0) dB
	+
	\frac{C_U^2}{4a^2 nq } \| X\|_F^4 \int \|(B-B_0) \|_F^4 		\pi(B - B_0) dB
	\\
	&	\leq 
	\frac{C_U}{a nq } \| X\|_F^2 qp \tau^2
	+
	\frac{C_U^2}{4a^2 nq } \| X\|_F^4 
	\left(
	\int \|(B-B_0) \|_F^2 \pi(B - B_0) dB
	\right)^2
	\\
	&	\leq 
	\frac{C_U}{a nq } \| X\|_F^2 qp \tau^2
	+
	\frac{C_U^2}{4a^2 nq } \| X\|_F^4 q^2 p^2 \tau^4	
	,
	\end{align*}
	where we have used Hölder's inequality and Lemma \ref{lemma:arnak:1}.
	From Lemma \ref{lemma:arnak:2}, we have that 
	\begin{align*}
	\frac{1}{nq} KL(\rho_n,\pi)
	\leq
	\frac{	2 r^* (q+p+2) \log \left( 1+ \frac{\| B_0 \|_F}{\tau \sqrt{2r^*}} \right) 
	}{nq}
	.
	\end{align*}
	For $\tau^2 = a/(qp \|X \|_{F}^2)$, it leads to that
	\begin{align*}
	\frac{1}{nq} KL(\rho_n,\pi)
	\leq
	\frac{	2 r^* (q+p+2) \log \left( 1+  \frac{\| X \|_F \| B_0 \|_F \sqrt{ qp } }{\sqrt{ 2 a r^*}} \right) 
	}{nq}
	,
	\end{align*}
	\begin{align*}
	\int KL (P_{B_0},P_{B})
	\rho_n(dB)
	\leq 
	\frac{C_U}{2nq} 
	,
	\end{align*}
	and
	\begin{align*}
	\int
	\mathbb{E}\left[\log
	\left(
	\frac{p_{B}}{p_{B_0}}(Y)
	\right)^2
	\right] 
	\rho_n(dB)
	\leq
	\frac{C_U}{nq} 
	+
	\frac{C_U^2}{4 nq} 
	.
	\end{align*}
	Putting
	\[
	\varepsilon_n
	=
	\frac{C_U^2}{4nq}
	\vee
	\frac{ C_U	2 r^* (q+p+2) \log \left( 1+ \frac{\| X \|_F \| B_0 \|_F  \sqrt{ qp}}{\sqrt{ 2a r^*}}  \right) 
	}{nq}
	,
	\]
	we can apply Theorem 2.4 and using Corollary 2.5 in \cite{alquier2020concentration} to obtain our result. 
	This completes the proof.
	
\end{proof}

\begin{proof}[\bf Proof of Proposition \ref{propos_estimation2}]
	From Theorem \ref{thm_result_probabily}, one has that
	\begin{equation*}
	\mathbb{P} \left[ \int D_{\alpha}(P_{B},P_{B_0}) \pi_{n,\alpha}({\rm d} B ) 
	\leq 
	2	\frac{1+\alpha}{1-\alpha}\varepsilon_n
	\right]
	\geq 
	1 -	\frac{2}{n\varepsilon_n}
	,
	\end{equation*}
	and from Lemma \ref{lema_lowerboundreymi}, where $ \theta = X^\top B , \theta_0 = X^\top B_0 $, one has that
	\begin{align*}
	D_{\alpha}(P_{B},P_{B_0})
	\geq
	\frac{C_L }{2a nq}  
	\|X^\top (B -B_0) \|_F^2
	.
	\end{align*}
	Thus, we obtain that
	\begin{equation*}
	\mathbb{P} \left[ \int 	
	\frac{C_L }{2a nq}  
	\| X^\top(B - B_0) \|_F^2 \pi_{n,\alpha}({\rm d} B ) 
	\leq 
	2	\frac{1+\alpha}{1-\alpha}\varepsilon_n
	\right]
	\geq 
	1 -	\frac{2}{n\varepsilon_n}
	,
	\end{equation*}	
	consequently,
	\begin{equation*}
	\mathbb{P} \left[ \int  
	\frac{1}{nq}
	\| X^\top(B - B_0)\|_F^2 \pi_{n,\alpha}({\rm d} B ) 
	\leq 
	\frac{ 4a (1+\alpha)}{ C_L (1- \alpha) }\varepsilon_n
	\right]
	\geq 
	1 -	\frac{2}{n\varepsilon_n}
	.
	\end{equation*}	
	Applying Jensen's inequality with a convex function, we find that
	\begin{equation*}
\mathbb{P} \left[  
	\frac{1}{nq} 
\| X^\top(\hat{B} - B_0)\|_F^2 
\leq 
\frac{ 4a (1+\alpha)}{ C_L (1-\alpha) }\varepsilon_n
\right]
\geq 
1 -	\frac{2}{n\varepsilon_n}
.
\end{equation*}		
	The proof is completed.	
	
\end{proof}

\subsection{Lemmas}

\begin{lemma}
	\label{lema_KLupperbound}
	Let Assumption \ref{assume_boundon2ndorder} hold. Then, for any $ \theta, \zeta \in \mathbb{R}^{n \times q} $, we have that
	\begin{equation*}
	KL(P_{\theta},P_{\zeta} )
	\leq
	\frac{C_U}{2a}  \frac{1}{nq} \| \zeta -\theta \|_F^2
	.
	\end{equation*}
\end{lemma}

\begin{proof}[\bf Proof of Lemma \ref{lema_KLupperbound}]
The rationale of the proof is based on Lemma 1 from \cite{abramovich2016model}. In our model \eqref{eq_main_model}, one has that
\begin{align}
	KL(P_{\theta},P_{\zeta} )
&	=
		\mathbb{E}\left[\log
\left(
\frac{p_{\theta}}{p_\zeta}(Y)
\right)
\right] 
\nonumber
\\
& = 
\frac{1}{a} \frac{1}{nq}
\sum_{ij} \left[\mathbb{E}
 \left(
Y_{ij} (\theta_{ij} - \zeta_{ij}) - b(\theta_{ij}) + b(\zeta_{ij})
\right)
\right] 
\nonumber
\\
& = 
\frac{1}{a}  \frac{1}{nq}
\sum_{ij}\left[
b'(\theta_{ij}) (\theta_{ij} - \zeta_{ij}) - b(\theta_{ij}) + b(\zeta_{ij})
\right]
\label{eq_kl01}
.
\end{align}
Utilizing a Taylor expansion of \( b(\zeta_{ij}) \) around \( \theta_{ij} \), it follows that
$$
b(\zeta_{ij}) 
= 
b(\theta_{ij}) + b'(\theta_{ij}) (\zeta_{ij} -\theta_{ij}) 
+ 
\frac{b''(c_{ij})}{2}  (\zeta_{ij} -\theta_{ij})^2
,
$$
where $ c_{ij} $ lies between $ \theta_{ij} $ and $ \zeta_{ij} $. Substituting this into \eqref{eq_kl01}, we have that
\begin{align*}
	KL(P_{\theta},P_{\zeta} )
	=
\frac{1}{2a}  \frac{1}{nq}
\sum_{ij} b''(c_{ij})  (\zeta_{ij} -\theta_{ij})^2
.
\end{align*}
Under the Assumption \ref{assume_boundon2ndorder}, we have that $ 	\sup_{s\in\mathbb{R}}  b''(s) \leq C_U $. Therefore, one obtains that
	\begin{equation*}
KL(P_{\theta},P_{\zeta} )
\leq
\frac{C_U}{2a}  \frac{1}{nq}
\sum_{ij} (\zeta_{ij} -\theta_{ij})^2
.
\end{equation*}
This completes the proof.

\end{proof}

\begin{lemma}
	\label{lema_KLup_missped}
	Let Assumption \ref{assume_bound1storder} hold. Then, for any $ \theta, \zeta \in \mathbb{R}^{n \times q} $, we have that
	\begin{equation*}
	\mathbb{E}_{ \theta_0}
	\left[\log
	\left(
	\frac{p_{\bar{\theta}}}{p_\zeta}(Y)
	\right)
	\right]
	\leq
	\frac{2U_1}{a} \frac{1}{ \sqrt{nq} }
   \| \zeta - \bar\theta \|_F
	.
	\end{equation*}
\end{lemma}
\begin{proof}[Proof of Lemma \ref{lema_KLup_missped}]
	One have that
	\begin{align}
	\mathbb{E}_{ \theta_0}
	\left[\log
	\left(
	\frac{p_{\bar{\theta}}}{p_\zeta}(Y)
	\right)
	\right]
	& = 
	\frac{1}{a} \frac{1}{nq}
	\sum_{ij}
	\left[ 	\mathbb{E}_{ \theta_0}
	\left(
	Y_{ij} (\bar{\theta}_{ij} - \zeta_{ij}) - b(\bar{\theta}_{ij}) + b(\zeta_{ij})
	\right)
	\right] 
	\nonumber
	\\
	& = 
	\frac{1}{a}\frac{1}{nq}
	\sum_{ij}\left[
	b'(\theta_{0,ij}) (\bar{\theta}_{ij} - \zeta_{ij}) - b(\bar{\theta}_{ij}) + b(\zeta_{ij})
	\right]
	\label{eq_KL_misped1}
	.
	\end{align}
	Utilizing a Taylor expansion of \( b(\zeta_{ij}) \) around \( \bar{\theta}_{ij} \), it follows that
	$$
	b(\zeta_{ij}) 
	= 
	b(\bar{\theta}_{ij}) + b'(c_{ij}) (\zeta_{ij} -\bar{\theta}_{ij}) 
	,
	$$
	where $ c_{ij} $ lies between $ \theta_{ij} $ and $ \zeta_{ij} $. Substituting this into \eqref{eq_KL_misped1}, we have that
	\begin{align*}
	\mathbb{E}_{ \theta_0}
	\left[\log
	\left(
	\frac{p_{\bar{\theta}}}{p_\zeta}(Y)
	\right)
	\right]
&	=
	\frac{1}{a}\frac{1}{nq}
	\sum_{ij}\left[
	b'(\theta_{0,ij}) (\bar{\theta}_{ij} - \zeta_{ij})
	+ b'(c_{ij}) (\zeta_{ij} -\bar{\theta}_{ij}) 
	\right]
	\\	
&	\leq
	\frac{1}{a}\frac{1}{nq}
	\sum_{ij}\left[
	| b'(\theta_{0,ij}) (\bar{\theta}_{ij} - \zeta_{ij}) |
	+
	| b'(c_{ij}) (\zeta_{ij} -\bar{\theta}_{ij}) |
	\right]
	.
	\end{align*}
	Under the Assumption \ref{assume_bound1storder}, we have that $ 	\sup_{s\in\mathbb{R}}  b'(s) \leq U_1 $. Therefore, one obtains that
	\begin{align*}
	\mathbb{E}_{ \theta_0}
	\left[\log
	\left(
	\frac{p_{\bar{\theta}}}{p_\zeta}(Y)
	\right)
	\right]
	&	\leq
	\frac{2U_1}{a}\frac{1}{nq}
	\sum_{ij} |\zeta_{ij} - \bar\theta_{ij}|
	\\
	& \leq
	\frac{2U_1}{a}\frac{1}{nq}
	\sqrt{nq}   \| \zeta - \bar\theta \|_F
	.
	\end{align*}
	This completes the proof.
	
\end{proof}

\begin{lemma}
	\label{lema_lowerboundreymi}
	Let Assumption \ref{assume_lowerbound2order} hold. Then, for any $ \theta, \zeta \in \mathbb{R}^{n \times q} $, we have that
	\begin{align*}
	D_{\alpha}(P_{\theta},P_{\zeta})
	\geq
	\frac{C_L }{2a }  \frac{1}{nq}
	\|\theta -\zeta \|_F^2
	.
	\end{align*}
\end{lemma}

\begin{proof}[\bf Proof of Lemma \ref{lema_lowerboundreymi}]
	As we can write that
	\begin{align*}
	D_{\alpha}(P_{\theta},P_{\zeta})	
=
	\frac{1}{\alpha-1} \log 
	\int \left(\frac{{\rm d} P_{\theta} }{{\rm d} P_{\zeta} }\right)^{\alpha -1 }
	{\rm d} P_{\theta}  
	\geq
 \int 
	\log \left(
	\frac{{\rm d} P_{\theta} }{{\rm d} P_{\zeta} }\right)
	{\rm d} P_{\theta}  
	\end{align*}
	where we have used Jensen's inequality and the concavity of the logarithmic function.
	Now, 
	\begin{align*}
	\int 
	\log \left(
	\frac{{\rm d} P_{\theta} }{{\rm d} P_{\zeta} }\right)
	{\rm d} P_{\theta}  
	&	=
	\frac{1}{a} \frac{1}{nq} 
	\sum_{ij} \int 
	\left(
	Y_{ij} (\theta_{ij} - \zeta_{ij}) - b(\theta_{ij}) + b(\zeta_{ij})
	\right)
	{\rm d} P_{\theta_{ij}}  
	\\
	& =	
	\frac{1}{a} \frac{1}{nq} \sum_{ij}	\left[
	b'(\theta_{ij})(\theta_{ij} - \zeta_{ij}) - b(\theta_{ij}) + b(\zeta_{ij})
	\right]
	\end{align*}
	A Taylor expansion of \( b(\zeta_{ij}) \) around \( \theta_{ij} \), it follows that
	$$
	b(\zeta_{ij}) 
	= 
	b(\theta_{ij}) + b'(\theta_{ij}) (\zeta_{ij} -\theta_{ij}) 
	+ 
	\frac{b''(c_{ij})}{2}  (\zeta_{ij} -\theta_{ij})^2
	,
	$$
	where $ c_{ij} $ lies between $ \theta_{ij} $ and $ \zeta_{ij} $. Thus we have
	\begin{align*}
	D_{\alpha}(P_{\theta},P_{\zeta})
	\geq
	\frac{1}{a} \frac{1}{nq} 
	\sum_{ij}
	\frac{b''(c_{ij})}{2}  (\theta_{ij} -\zeta_{ij})^2.
	\end{align*}	
Using Assumption \ref{assume_lowerbound2order}, we have that
	\begin{align*}
	D_{\alpha}(P_{\theta},P_{\zeta})
	\geq
	\frac{C_L }{2a }  \frac{1}{nq}
	\|\theta -\zeta\|_F^2
	.
	\end{align*}
	This completes the proof.
\end{proof}

\begin{definition}
	\label{dfn:posterior:transla}
	Let's define 
	\begin{equation}
	\label{eq_priorspecific_anark}
	\rho_{0}(B) = \pi(B - B_0)
	.	
	\end{equation}
\end{definition}

\begin{lemma}[Lemma 1 in \cite{dalalyan2020exponential}]
	\label{lemma:arnak:1}
	We have
	$$
	\int \| B \|_{F}^2 \pi( B ){\rm d} B  \leq q p \tau^2. 
	$$
\end{lemma}

\begin{lemma}[Lemma 2 in \cite{dalalyan2020exponential}]
	\label{lemma:arnak:2}
	We have
	$$
	KL( \rho_{0 } ,  \pi) 
	\leq 
	2 {\rm rank} (B_0 ) (q+p+2) \log \left( 1+ \frac{\| B_0 \|_F}{\tau \sqrt{2{\rm rank} (B_0 )}} \right) 
	,
	$$
	with the convention $0\log(1+0/0)=0$.
\end{lemma}

\begin{lemma}
	\label{lm_bound_varkl}
	Assume that Assumption \ref{assume_boundon2ndorder} holds. Then, for any $ \theta, \zeta \in \mathbb{R}^{n \times q} $, we have that
	\begin{align*}
	\mathbb{E}\left[\log
	\left(
	\frac{p_{\theta}}{p_\zeta}(Y)
	\right)^2
	\right]
	\leq 
	\frac{C_U}{a} \frac{1}{nq}  \| \theta - \zeta \|_F^2
	+
	\frac{C_U^2}{4a^2} \frac{1}{nq}  \| \zeta -\theta\|_F^4
	.
	\end{align*}
\end{lemma}

\begin{proof}[{\bf Proof of Lemma \ref{lm_bound_varkl}}]
	In our GLM in \eqref{eq_main_model}, one has that
	\begin{align}
	&	\mathbb{E}\left[\log
	\left(
	\frac{p_{\theta}}{p_\zeta}(Y)
	\right)^2
	\right] \nonumber
	\\
	& = 
	\frac{1}{a^2}  \frac{1}{nq} 
\sum_{ij}	
	\left[\mathbb{E} \left(
	Y_{ij} (\theta_{ij} - \zeta_{ij}) - b(\theta_{ij}) + b(\zeta_{ij})
	\right)^2\right] \nonumber
	\\
	& = 
	\frac{1}{a^2} \frac{1}{nq} 
	\sum_{ij}\mathbb{E}
	\left\{
	Y_{ij}^2 (\theta_{ij} - \zeta_{ij})^2
	-
	2 Y_{ij} (\theta_{ij} - \zeta_{ij}) [b(\theta_{ij}) - b(\zeta_{ij}) ]
	+
	[b(\theta_{ij}) - b(\zeta_{ij})]^2
	\right\}
	.	
	\label{eq_varkl_1}
	\end{align}
	We remind that from the GLM in \eqref{eq_main_model}, one has that $ \mathbb{E} (Y_{ij}) = b'(\theta_{ij}) $ and $ Var(Y_{ij}) = ab''(\theta_{ij}) $, see \cite{mccullagh1989generalized}. It leads to that 
	\begin{equation}
	\label{eq_ey2}
	\mathbb{E} Y_{ij}^2
	= 
	[\mathbb{E} (Y_{ij})]^2 +  Var(Y_{ij}) 
	= 
	[b'(\theta_{ij})]^2 + ab''(\theta_{ij}) 
	.
	\end{equation}
	Plugging \eqref{eq_ey2} into \eqref{eq_varkl_1}, one gets that
	\begin{align}
	& a^2 nq	
	\mathbb{E}\left[\log
	\left(
	\frac{p_{\theta}}{p_\zeta}(Y)
	\right)^2
	\right] \nonumber
	\\
	 =  &
	\sum_{ij} \left\{
	[b'(\theta_{ij})]^2 + ab''(\theta_{ij}) 
	\right\} (\theta_{ij} - \zeta_{ij})^2
 \nonumber
\\
  &
\hspace*{.5cm}	-
	\sum_{ij} 	2 b'(\theta_{ij})  (\theta_{ij} - \zeta_{ij}) [b(\theta_{ij}) - b(\zeta_{ij}) ]
	+
\sum_{ij}	[b(\theta_{ij}) - b(\zeta_{ij})]^2
	\nonumber
	\\
	 =&  
	\sum_{ij}
	ab''(\theta_{ij}) 
	(\theta_{ij} - \zeta_{ij})^2
	+ \sum_{ij}
	\left[ b'(\theta_{ij})(\theta_{ij} - \zeta_{ij}) - b(\theta_{ij}) + b(\zeta_{ij}) \right]^2
	\label{eq_varkl_2}
	.	
	\end{align}
	Using a Taylor expansion of $ b(\zeta_{ij}) $ around $ \theta_{ij} $, one has that 
	$$
	b(\zeta_{ij}) 
	= 
	b(\theta_{ij}) + b'(\theta_{ij}) (\zeta_{ij} -\theta_{ij}) 
	+ 
	\frac{b''(c_{ij})}{2}  (\zeta_{ij} -\theta_{ij})^2
	,
	$$
	where $ c_{ij} $ lies between $ \theta_{ij} $ and $ \zeta_{ij} $. Substituting this into \eqref{eq_varkl_2}, one gets that
	\begin{align*}
	a^2	nq
	\mathbb{E}\left[\log
	\left(
	\frac{p_{\theta}}{p_\zeta}(Y)
	\right)^2
	\right] 
	=  
	\sum_{ij}
	ab''(\theta_{ij}) 
	(\theta_{ij} - \zeta_{ij})^2
	+
	\sum_{ij}
	\frac{b''(c_{ij})^2}{2^2}  (\zeta_{ij} -\theta_{ij})^4
	.	
	\end{align*}
	By noting that $ \sum_{ij}
	(\zeta_{ij} -\theta_{ij})^4
	\leq 
	\left[ \sum_{ij}
	(\zeta_{ij} -\theta_{ij})^2 \right]^2
	$,
	we have, under Assumption \ref{assume_boundon2ndorder}, that
	\begin{align*}
	\mathbb{E}\left[\log
	\left(
	\frac{p_{\theta}}{p_\zeta}(Y)
	\right)^2
	\right] 
	\leq
	\frac{C_U}{a}  \frac{1}{nq}  \| \theta - \zeta \|_F^2
	+
	\frac{C_U^2}{4a^2}  \frac{1}{nq} \| \zeta -\theta\|_F^4
	.	
	\end{align*}  
	This completes the proof.
	
\end{proof}

\section{Closing Discussion}
\label{sc_discusconclue}

In this work, we provided consistency outcomes and convergence rates for the fractional posterior in generalized reduced rank regression by employing a low-rank promoting prior on the matrix parameter. We also derived corresponding results for the matrix parameter in terms of the Frobenius norm. Furthermore, we tackled model mis-specification, confirming the effectiveness of both the fractional posterior and its mean estimator. More specifically, we introduced an oracle-type inequality for generalized reduced rank regression.

While our current investigation primarily delves into theoretical aspects, we now provide a brief commentary on the computation of the fractional posterior. Given that the spectral scaled Student prior in \eqref{prior_scaled_Student} operates on the entire parameter matrix, one can utilize Langevin Monte Carlo methods to sample from the fractional posterior, as outlined in \cite{dalalyan2020exponential}.

Additionally, it is worth noting that by employing another type of prior known as low-rank factorization priors, which are popular in low-rank matrix estimation problems, similar concentration results for the fractional posterior can be achieved, albeit with more stringent assumptions \citep{alquier2020concentration,mai2024concentration}. However, since a prior is typically not conjugate in generalized linear models, one may need to explore MCMC methods when using low-rank factorization priors. The development of fast approximation methods based on optimization such as Variational Inference for low-rank factorization priors in our reduced-rank generalized linear models would present an intriguing avenue for future research.

\subsection*{Acknowledgments}
The author acknowledge support from the Norwegian Research Council, grant number 309960, through the Centre for Geophysical Forecasting at NTNU. The author expresses gratitude to Pierre Alquier for his valuable discussions regarding his paper on AoS.

\subsection*{Conflicts of interest/Competing interests}
The author declares no potential conflict of interests.


\end{document}